# Domino Tilings of Graph Grids

VALCHO MILCHEV

This article is dedicated to domino tilings of certain types of graph grids. For each of these grids, the domino tilings are represented using linear-recurrent sequences. New dependencies are proved that are not included in Neil Sloane's Online Encyclopedia of Integer Sequences (OEIS).

**1. Introduction**

We will note three linear-recurrent sequences, which have presence in the counting of the domino tilings: Fibonacci numbers, Pell numbers and Heronian triples.

We will denote the Fibonacci numbers with $f_n$. They are defined with $f_0 = 1$, $f_1 = 1$ and $f_n = f_{n-1} + f_{n-2}$ when $n \geq 2$. In OEIS the Fibonacci sequence is with number A000045. The domino tilings of a square grid $2 \times n$ and respectively the domino tilings of a graph grid $P_2 \times P_n$, $n \geq 1$, shown on Figure 1, are expressed by Fibonacci numbers $\{f_n\} = \{1, 2, 3, 5, 8, ...\}$:

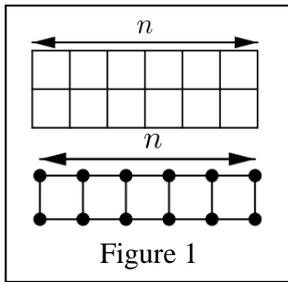

Figure 1

$$f_n = \frac{1}{\sqrt{5}}\left[\left(\frac{1+\sqrt{5}}{2}\right)^{n+1} + \left(\frac{1-\sqrt{5}}{2}\right)^{n+1}\right].$$

We will note the Pell numbers with $p_n$. They are defined with $p_0 = 1$, $p_1 = 2$ and $p_n = 2p_{n-1} + p_{n-2}$ when $n \geq 2$. They are solutions for $x$ of Pell's equation $x^2 - 2y^2 = 1$. The sequence of Pell's numbers $\{p_n\} = \{1, 2, 5, 12, 29, 70, 169, ...\}$ in OEIS is with number A000129. The corresponding solutions for $y$ of Pell's equation $x^2 - 2y^2 = 1$ represent the sequence of numbers $\{h_n\} = \{1, 1, 3, 7, 17, 41, 99, ...\}$, where $h_0 = 1$, $h_1 = 1$ and $h_n = 2h_{n-1} + h_{n-2}$ when $n \geq 2$. In OEIS this sequence is with number A001333.

One more number sequence has attracted the mathematicians' attention for many centuries. We will call them the Heronian numbers. These are the numbers $\{H_n\} = \{2, 4, 14, 52, ...\}$, where $H_0 = 2$, $H_1 = 4$, $H_n = 4H_{n-1} - H_{n-2}$, when $n \geq 2$. They are remarkable with this that a triangle with sides $H_n - 1$, $H_n$ and $H_n + 1$ has an area that is an integer when $n \geq 1$. The triple $(H_n - 1; H_n; H_n + 1)$ is also called a Heronian triple, $n \geq 1$. The numbers $H_n$ are solution for $y$ of Pell's equation $x^2 - 3y^2 = 3$. In OEIS this sequence is with number A003500.

## 2. Ribbon Grids and Heronian Triples

We will consider a domino tiling of the graph grids, shown on Figure 2. Let's denote the graph grids and the number of the domino tilings for the same grids with $A_n$, $B_n$, $C_n$, $D_n$ and $E_n$.

**Theorem 1.** Let's consider the ribbon graph grids, shown on Figure 2. We have the following dependencies between them:

$A_n = B_n + C_n$;
$B_n = A_{n-1} + B_{n-1}$;
$C_n = A_{n-1} + B_n$;
$A_n = A_{n-1} + 2B_n$;
$A_{n-1} + A_n = 2C_n$;
$D_n = B_n + D_{n-1}$;
$E_n = B_n + E_{n-1}$.

**Proof.** The first equation is represented on Figure 3. For the top left vertex there are exactly two possibilities for the connected edge – highlighted in red. After specifying the first edge, the edge in green is marked as "obligatory". For the residual graph, the domino tilings are $B_n$. In the other case for the first marked edge for the residual graph, the domino tilings are $C_n$. In a similar way to figures 4, 5, 6 and 7 are illustrated more dependencies.

Using the proved equations in theorem 1 we compose table 1 for the domino tilings.

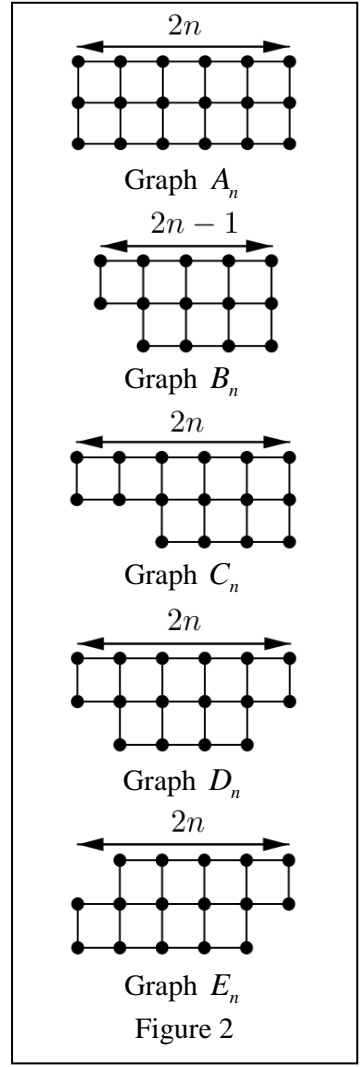

Graph $A_n$

Graph $B_n$

Graph $C_n$

Graph $D_n$

Graph $E_n$

Figure 2

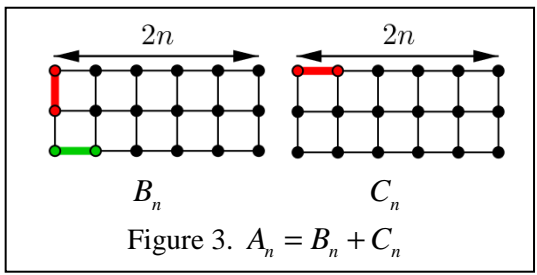

Figure 3. $A_n = B_n + C_n$

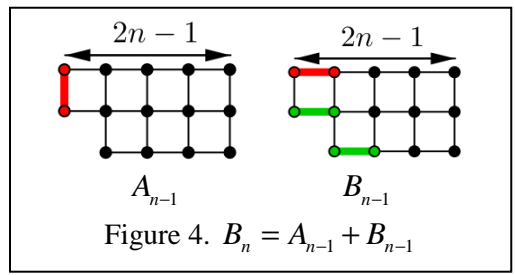

Figure 4. $B_n = A_{n-1} + B_{n-1}$

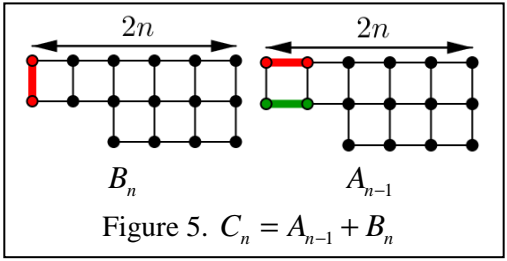

Figure 5. $C_n = A_{n-1} + B_n$

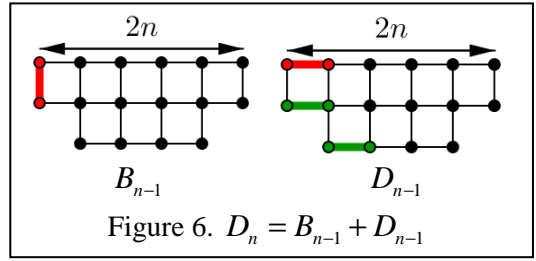

Figure 6. $D_n = B_{n-1} + D_{n-1}$

| $n$ | $A_n$ | $B_n$ | $C_n$ | $D_n$ | $E_n$ |
|---|---|---|---|---|---|
| 1 | 3 | 1 | 2 | 2 | 1 |
| 2 | 11 | 4 | 7 | 6 | 5 |
| 3 | 41 | 15 | 26 | 21 | 20 |
| 4 | 153 | 56 | 121 | 77 | 76 |

Table 1

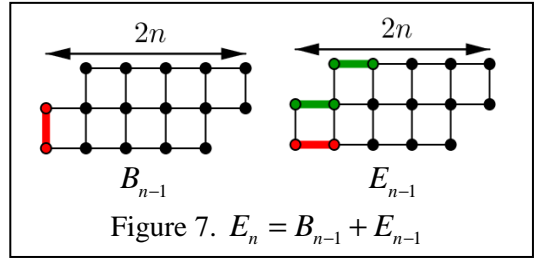

Figure 7. $E_n = B_{n-1} + E_{n-1}$

**Theorem 2.** The number of the domino tilings for the ribbon grids shown on figure 2 is represented for $n \geq 1$, using the formulas:

$$A_n = \frac{1}{6}\left[(3+\sqrt{3})(2+\sqrt{3})^n + (3-\sqrt{3})(2-\sqrt{3})^n\right];$$

$$B_n = \frac{1}{2\sqrt{3}}\left[(2+\sqrt{3})^n - (2-\sqrt{3})^n\right];$$

$$C_n = \frac{1}{2}\left[(2+\sqrt{3})^n + (2-\sqrt{3})^n\right];$$

$$D_n = \frac{1}{12}\left[(3+\sqrt{3})(2+\sqrt{3})^n + (3-\sqrt{3})(2-\sqrt{3})^n\right] + \frac{1}{2};$$

$$E_n = \frac{1}{12}\left[(3+\sqrt{3})(2+\sqrt{3})^n + (3-\sqrt{3})(2-\sqrt{3})^n\right] - \frac{1}{2}.$$

**Proof.** Using equations $B_n = A_{n-1} + B_{n-1}$ and $A_n = A_{n-1} + 2B_n$ we obtain $A_n = 4A_{n-1} - A_{n-2}$, $B_n = 4B_{n-1} - B_{n-2}$, $n \geq 3$. The sequences $\{A_n\}$ and $\{B_n\}$ are linear-recurrent and they have the same characteristic equation $x^2 - 4x + 1 = 0$ with roots $x_1 = 2+\sqrt{3}$ and $x_2 = 2-\sqrt{3}$. We look for these kind of formulas $A_n = c_1 x_1^n + c_2 x_2^n$. We determine the coefficients using the data from table 1. Hence, for $n \geq 1$ we have:

$$A_n = \frac{1}{6}\left[(3+\sqrt{3})(2+\sqrt{3})^n + (3-\sqrt{3})(2-\sqrt{3})^n\right],$$

$$B_n = \frac{1}{2\sqrt{3}}\left[(2+\sqrt{3})^n - (2-\sqrt{3})^n\right].$$

From the equation $C_n = A_{n-1} + B_n$ we receive

$$C_n = 4C_{n-1} - C_{n-2} \text{ за } n \geq 3 \text{ and } C_n = \frac{1}{2}\left[(2+\sqrt{3})^n + (2-\sqrt{3})^n\right] \text{ for } n \geq 1.$$

We receive the following equations for sequence $\{D_n\}$, when $n \geq 4$

$$D_n - D_{n-1} = B_n = 4B_{n-1} - B_{n-2} = 4(D_{n-1} - D_{n-2}) - (D_{n-2} - D_{n-3}),$$

$$D_n = 5D_{n-1} - 5D_{n-2} + D_{n-3}.$$

By analogy, for sequence $\{E_n\}$, we have $E_n = 5E_{n-1} - 5E_{n-2} + E_{n-3}$. Using table 1 we prove the equations:

$$D_n = 4D_{n-1} - D_{n-2} - 1 \text{ and } E_n = 4E_{n-1} - E_{n-2} + 1, \ n \geq 3.$$

The characteristic equation of sequences $\{D_n\}$ and $\{E_n\}$ is $x^3 - 5x^2 + 5x - 1 = 0$, the roots of which are $x_1 = 2+\sqrt{3}$, $x_2 = 2-\sqrt{3}$ and $x_3 = 1$. Using the standard way for the formula for third-order linear-recurrent sequence we receive

$$D_n = \frac{1}{12}\left[(3+\sqrt{3})(2+\sqrt{3})^n + (3-\sqrt{3})(2-\sqrt{3})^n\right] + \frac{1}{2}.$$

$$E_n = \frac{1}{12}\left[(3+\sqrt{3})(2+\sqrt{3})^n + (3-\sqrt{3})(2-\sqrt{3})^n\right] - \frac{1}{2}.$$

**Corollary.** The proved formulas give the opportunity for receiving new dependencies:

$$D_n = E_n + 1; \ 2D_n = A_n + 1; \ 2E_n = A_n - 1; \ A_n = D_n + E_n.$$

**Note.** We notice that $2C_n = H_n$ - the middle-sized side of the Heronian triple. This means that $A_n + A_{n-1} = H_n$,

### 3. Domino tilings of graph grids $C_4 \times P_n$

**Theorem 3.** Let $G_n$ be the number of the domino tilings of the graph $C_4 \times P_n$, shown on figure 8. Let $A_n$ and $B_n$ be the numbers from theorem 2. Then

$$G_n = \frac{1}{6}\left[(2+\sqrt{3})^{n+1} + (2-\sqrt{3})^{n+1}\right] + \frac{1}{3}(-1)^n;$$

$$G_{2n} = A_n^2, \ G_{2n-1} = 2B_n^2, \ n \geq 1.$$

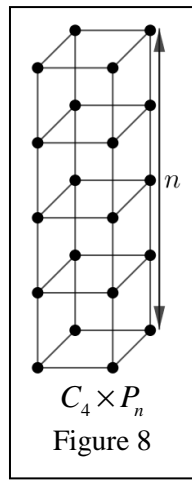

$C_4 \times P_n$
Figure 8

**Proof.** Let's note with $g_n$ the graph grid on figure 9, and also the number of the domino tilings for the same graph. Then $G_n = 2G_{n-1} + G_{n-2} + 4g_{n-1}$. This recurrent equation is illustrated on figure 10. All the ways for the "first" highlighted edges are shown – two horizontal, four vertical and a combination of two horizontal and two vertical.

From the other hand, we have $g_n = G_{n-1} + g_{n-1}$ - this equation is illustrated on figure 11. Using the received two recurrent equations we find

$$G_n = 3G_{n-1} + 3G_{n-2} - G_{n-3}, \quad g_n = 3g_{n-1} + 3g_{n-2} - g_{n-3}, \quad n \geq 4.$$

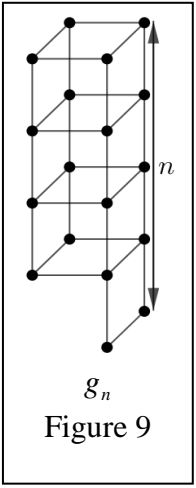

$g_n$

Figure 9

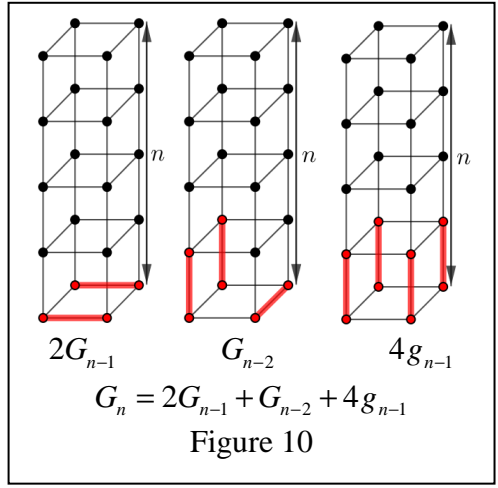

$2G_{n-1} \qquad G_{n-2} \qquad 4g_{n-1}$

$G_n = 2G_{n-1} + G_{n-2} + 4g_{n-1}$

Figure 10

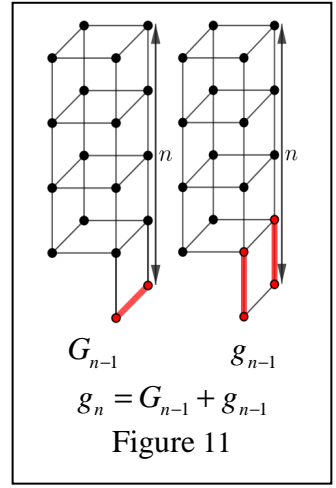

$G_{n-1} \qquad g_{n-1}$

$g_n = G_{n-1} + g_{n-1}$

Figure 11

The sequences $\{G_n\}$ and $\{g_n\}$ have the same characteristic equation $x^3 - 3x^2 - 3x + 1 = 0$. Its roots are $x_1 = 2 + \sqrt{3}$, $x_2 = 2 - \sqrt{3}$ and $x_3 = -1$. Using that $G_1 = 2$, $G_2 = 9$, $G_3 = 32$, $g_1 = 1$, $g_2 = 3$, $g_3 = 12$ we get the following formulas for $n \geq 1$

$$G_n = \frac{1}{6}\left(2+\sqrt{3}\right)^{n+1} + \frac{1}{6}\left(2-\sqrt{3}\right)^{n+1} + \frac{1}{3}(-1)^n;$$

$$g_n = \frac{1}{12}\left[\left(1+\sqrt{3}\right)\left(2+\sqrt{3}\right)^n + \left(1-\sqrt{3}\right)\left(2-\sqrt{3}\right)^n\right] - \frac{1}{6}(-1)^n.$$

Now we easily prove the following:

$$A_n^2 = \left\{\frac{1}{6}\left[\left(3+\sqrt{3}\right)\left(2+\sqrt{3}\right)^n + \left(3-\sqrt{3}\right)\left(2-\sqrt{3}\right)^n\right]\right\}^2 =$$

$$= \frac{1}{6}\left(2+\sqrt{3}\right)^{2n+1} + \frac{1}{6}\left(2-\sqrt{3}\right)^{2n+1} + \frac{1}{3} = G_{2n},$$

$$2B_n^2 = 2\left\{\frac{1}{2\sqrt{3}}\left[\left(2+\sqrt{3}\right)^n - \left(2-\sqrt{3}\right)^n\right]\right\}^2 =$$

$$= \frac{1}{6}\left(2+\sqrt{3}\right)^{2n} + \frac{1}{6}\left(2-\sqrt{3}\right)^{2n} - \frac{1}{3} = G_{2n-1}.$$

**Corollary.** Let $A_n$, $B_n$, $G_n$ and $g_n$ be the numbers from theorem 1, theorem 2 and theorem 3. Then:

$3G_{2n} = A_{2n+2} + B_{2n+2} + 1$;

$3G_{2n-1} = A_{2n} + B_{2n} - 1$;

$A_n = G_n + G_{n-1} = g_{n+1} - g_{n-1}$.

| $n$ | $A_n$ | $B_n$ | $G_n$ | $g_n$ |
|---|---|---|---|---|
| 1 | 3 | 1 | 2 | 1 |
| 2 | 11 | 4 | 9 | 3 |
| 3 | 41 | 15 | 32 | 12 |
| 4 | 153 | 56 | 121 | 44 |
| 5 | 571 | 209 | 450 | 165 |
| 6 | 2 131 | 780 | 1 681 | 615 |
| 7 | 7 953 | 2 911 | 6 272 | 2 296 |
| 8 | 29 681 | 10 864 | 23 409 | 8 568 |
| 9 | 110 771 | 40 545 | 87 362 | 31 977 |
| 10 | 413 403 | 151 316 | 326 041 | 119 339 |

Table 2

In table 2 are shown the domino tilings $A_n$, $B_n$, $G_n$ and $g_n$ for the values of $n$ ranging from 1 to 10. The proved dependencies between them can be seen.

| $k$ | $A_k$ | $B_k$ | $k$ |  | $g_{2m+1}$ |  | $g_{2m}$ | $n$ |
|---|---|---|---|---|---|---|---|---|
|  |  |  | 1 |  | 1 |  |  | 1 |
| 1 | 3 | 1 | 1 |  |  | 3.1 | 3 | 2 |
| 1 | 3 | 4 | 2 | 3.4 | 12 |  |  | 3 |
| 2 | 11 | 4 | 2 |  |  | 11.4 | 44 | 4 |
| 2 | 11 | 15 | 3 | 11.15 | 165 |  |  | 5 |
| 3 | 41 | 15 | 3 |  |  | 41.15 | 615 | 6 |
| 3 | 41 | 56 | 4 | 41.56 | 2 296 |  |  | 7 |
| 4 | 153 | 56 | 4 |  |  | 153.56 | 8 568 | 8 |
| 4 | 153 | 209 | 5 | 153.209 | 31 977 |  |  | 9 |
| 5 | 571 | 209 | 5 |  |  | 571.209 | 119 339 | 10 |
| 5 | 571 | 780 | 6 | 571.780 | 445 380 |  |  | 11 |
| 6 | 2131 | 780 | 6 |  |  | 2131.780 | 1 662 180 | 12 |

Table 3

*Table 3 illustrates an interesting relation of $g_n$ with $A_n$ and $B_n$. The numbers $g_n$ are received when we multiply the numbers $A_n$ and $B_n$ in a defined order.*

**Theorem 4.** We define a sequence of numbers $z_n$ as follows: $z_{2k-1} = A_k$ and $z_{2k} = B_k$. Then we have $g_1 = 1$, $g_n = z_n \cdot z_{n+1}$ for $n \geq 2$.

**Proof.** We have to prove that $g_n = A_k \cdot B_k$ when $n$ is an even number and $g_n = B_k \cdot A_{k+1}$ when $n$ is an odd number. Actually for the respective products we have:

$$A_k \cdot B_k = \frac{1}{6}\left[\left(3+\sqrt{3}\right)\left(2+\sqrt{3}\right)^k + \left(3-\sqrt{3}\right)\left(2-\sqrt{3}\right)^k\right].$$

$$\cdot \frac{1}{2\sqrt{3}}\left[\left(2+\sqrt{3}\right)^{k}-\left(2-\sqrt{3}\right)^{k}\right]=$$

$$=\frac{1}{12}\left[\left(1+\sqrt{3}\right)\left(2+\sqrt{3}\right)^{2k}+\left(1-\sqrt{3}\right)\left(2-\sqrt{3}\right)^{2k}\right]-\frac{1}{6}=g_{2k};$$

$$B_{k+1}\cdot A_{k}=\frac{1}{2\sqrt{3}}\left[\left(2+\sqrt{3}\right)^{k+1}-\left(2-\sqrt{3}\right)^{k+1}\right]\cdot$$

$$\cdot\frac{1}{6}\left[\left(3+\sqrt{3}\right)\left(2+\sqrt{3}\right)^{k}+\left(3-\sqrt{3}\right)\left(2-\sqrt{3}\right)^{k}\right]=$$

$$=\frac{1}{12}\left[\left(1+\sqrt{3}\right)\left(2+\sqrt{3}\right)^{2k+1}+\left(1-\sqrt{3}\right)\left(2-\sqrt{3}\right)^{2k+1}\right]+\frac{1}{6}=g_{2k+1}.$$

| Theorem 1, theorem 2 and theorem 3 | | Sequences in OEIS | Number |
|---|---|---|---|
| $A_n$ | 3, 11, 41, 153, 571, 2131,... $A_n = 4A_{n-1} - A_{n-2}$, | 1, 1, 3, 11, 41, 153, 571, 2131, 7953, 29681,... | A001835 |
| $B_n$ | 1, 4, 15, 56, 209, 780, 2911,... $B_n = 4B_{n-1} - B_{n-2}$ | 0, 1, 4, 15, 56, 209, 780, 2911, 10864, 40545,... | A001353 |
| $C_n$ | 2, 7, 26, 97, 362, 1351, 5042,... $C_n = 4C_{n-1} - C_{n-2}$ | 1, 2, 7, 26, 97, 362, 1351, 5042, 18817, 70226,... | A001075 |
| $D_n$ | 2, 6, 21, 77, 286, 1066, 3977,... $D_n = 5D_{n-1} - 5D_{n-2} + D_{n-3}$ | 1, 2, 6, 21, 77, 286, 1066, 3977, 14841, 55386,... | A101265 |
| $E_n$ | 1, 5, 20, 76, 285, 1065, 3976,... $E_n = 5E_{n-1} - 5E_{n-2} + E_{n-3}$ | 0, 1, 5, 20, 76, 285, 1065, 3976, 14840, 55385,... | A061278 |
| $G_n$ | 2, 9, 32, 121, 450, 1681, ,... $G_n = 3G_{n-1} + 3G_{n-2} - G_{n-3}$ | 1, 2, 9, 32, 121, 450, 1681, 6272, 23409, 87362,... | A006253 |
| $g_n$ | 1, 3, 12, 44, 165, 615, 2296,... $g_n = 3g_{n-1} + 3g_{n-2} - g_{n-3}$ | 0, 1, 3, 12, 44, 165, 615, 2 296, 8568, 31977, ... | A109437 |
| Table 4 | | | |

Table 4 references a part of the sequences from theorems 1, 2 and 3 and the corresponding denotations in OEIS.

**Note**. The relation of the domino tilings of a graph $g_n$ and the domino tilings $A_n$ and $B_n$ is not indicated in OEIS.

### 4. Domino tilings of graph grids $W_4 \times P_n$

Let's consider the graph grid $(K_4 - e) \times P_n$, also known as $W_4 \times P_n$, shown on figure 12. We will use the shorter denotation $W_4 \times P_n$.

**Theorem 5.** For the domino tilings $V_n$ of a graph $W_4 \times P_n$ we have the recurrent dependency $V_n = 2V_{n-1} + 7V_{n-2} + 2V_{n-3} - V_{n-4}$.

**Proof.** Let's consider two more graphs, shown on figure 13 – these are variants, in which two vertices are removed from the graph $W_4 \times P_n$. The non-zero first highlights of graph $W_4 \times P_n$ are shown on figure 14. The next recurrent relations are illustrated

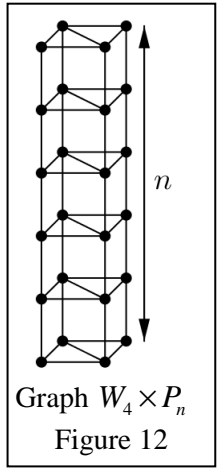

Graph $W_4 \times P_n$
Figure 12

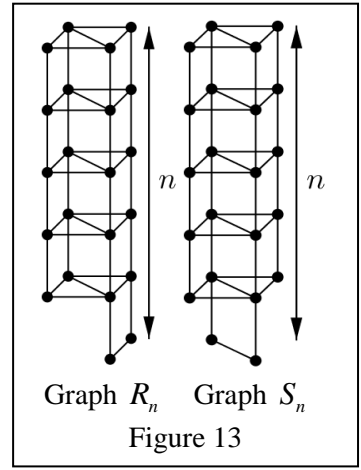

Graph $R_n$   Graph $S_n$
Figure 13

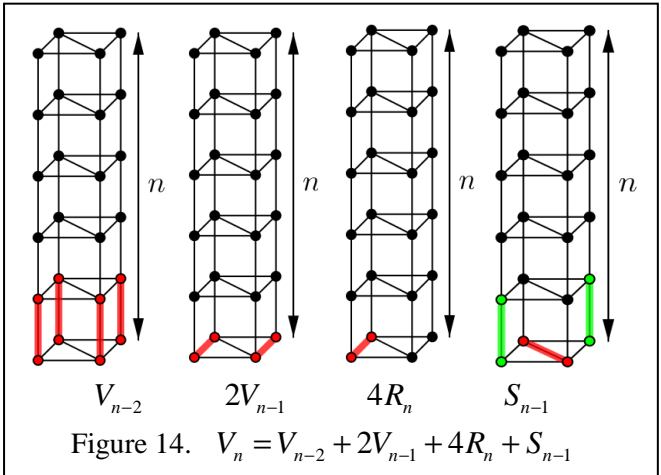

$V_{n-2}$  $2V_{n-1}$  $4R_n$  $S_{n-1}$
Figure 14.  $V_n = V_{n-2} + 2V_{n-1} + 4R_n + S_{n-1}$

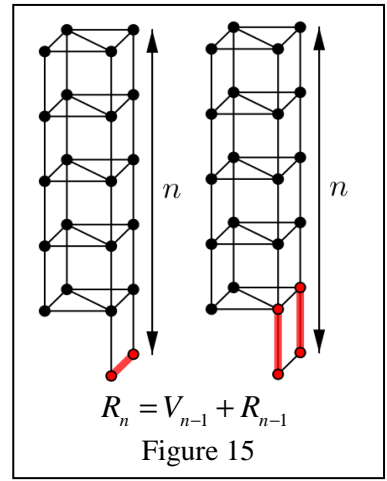

$R_n = V_{n-1} + R_{n-1}$
Figure 15

on figure 15 and figure 16. In this way we have a system of three recurrent equations:

$V_n = V_{n-2} + 2V_{n-1} + 4R_n + S_{n-1}$,

$R_n = V_{n-1} + R_{n-1}$,

$S_n = V_{n-1} + S_{n-2}$.

Then we obtain

$S_n - S_{n-2} = V_{n-1}$,

$R_n - R_{n-1} = V_{n-1}$,

$R_{n+1} - R_{n-1} = V_n + V_{n-1}$;

$V_n = V_{n-2} + 2V_{n-1} + 4R_n + S_{n-1}$;

$V_{n-2} = V_{n-4} + 2V_{n-3} + 4R_{n-4} + S_{n-3}$;

$V_n = 2V_{n-1} + 7V_{n-2} + 2V_{n-3} - V_{n-4}$.

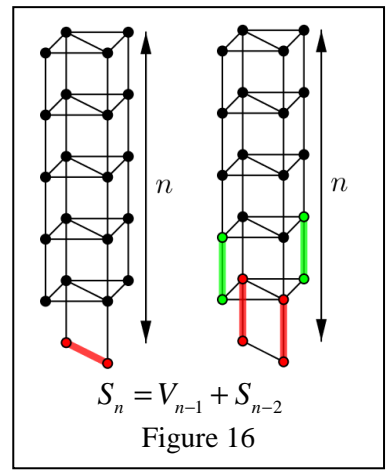

$S_n = V_{n-1} + S_{n-2}$
Figure 16

Using these equations we find $V_1 = 2$, $V_2 = 10$, $V_3 = 36$, $V_4 = 145$. We can find the explicit formula for $V_n$ using the solutions of the equation $x^4 - 2x^3 - 7x^2 - 2x + 1 = 0$. Its roots are

$$x_1 = \frac{(1+\sqrt{5})(1+\sqrt{2})}{2}, \quad x_2 = \frac{(1-\sqrt{5})(1-\sqrt{2})}{2},$$

$$x_3 = \frac{(1-\sqrt{5})(1+\sqrt{2})}{2}, \quad x_4 = \frac{(1+\sqrt{5})(1-\sqrt{2})}{2}.$$

The coefficients $c_1$, $c_2$, $c_3$, and $c_4$ of formula $V_n = c_1 x_1^n + c_2 x_2^n + c_3 x_3^n + c_4 x_4^n$, are defined by the initial conditions $V_1 = 2$, $V_2 = 10$, $V_3 = 36$, $V_4 = 145$. Hence, for $n \geq 1$ we will receive

$$V_n = \frac{\sqrt{10}}{40}\left(\left(\frac{1+\sqrt{5}}{2}\right)^n - \left(\frac{1-\sqrt{5}}{2}\right)^n\right)\left((1+\sqrt{2})^n - (1-\sqrt{2})^n\right).$$

However, considerable efforts are necessary in order to do this.

In 2002 James A. Sellers published one extraordinary formula for the domino tilings of a graph grid $W_4 \times P_n$, or namely: the number of the domino tilings of a graph $W_4 \times P_n$ is represented as a product of the respective numbers from the Fibonacci and Pell sequences. For some of the values this fact is illustrated in table 5. The proof is made using the method of mathematical induction.

**Theorem 6.** For the number of the domino tilings of a graph grid $W_4 \times P_n$ we have the formula $V_n = f_n \cdot p_n$, for $n \geq 1$.

**Note.** The formula $V_n = f_n \cdot p_n$ is a particular case of a more general statement: Let the linear-recurrent sequences $\{a_n\}$ and $\{b_n\}$ with recurrent equations $a_n = u \cdot a_{n-1} + a_{n-2}$ and $b_n = v \cdot b_{n-1} + b_{n-2}$ be given. Then the sequence with a common term $c_n = a_n \cdot b_n$ has the recurrent equation

$$c_n = uvc_{n-1} + (u^2 + v^2 + 2)c_{n-2} + uvc_{n-3} - c_{n-4}.$$

It follows that if $c_1 = a_1 \cdot b_1$, $c_2 = a_2 \cdot b_2$,

| $n$ | $f_n$ | $p_n$ | $V_n = f_n \cdot p_n$ |
|---|---|---|---|
| 1 | 1 | 2 | 2 |
| 2 | 2 | 5 | 10 |
| 3 | 3 | 12 | 36 |
| 4 | 5 | 29 | 145 |
| 5 | 8 | 70 | 560 |
| 6 | 13 | 169 | 2197 |
| 7 | 21 | 408 | 8568 |
| 8 | 34 | 985 | 33490 |

Table 5

$c_3 = a_3 \cdot b_3$, $c_4 = a_4 \cdot b_4$, then the sequence $\{c_n\}$ is with a recurrent equation

$$c_n = uvc_{n-1} + (u^2 + v^2 + 2)c_{n-2} + uvc_{n-3} - c_{n-4}.$$

In our case $u = 1$ and $v = 2$.

*As there is a graph (this is a graph grid $2 \times n$), for which the domino tilings are represented using the respective numbers from the Fibonacci sequence, the*

*question arises whether there is a graph grid, for which the domino tilings are represented by Pell's numbers. We will propound such a graph.*

**Theorem 7.** Let's consider two graph grids, shown on figure 17. We have the following dependencies:

$$p_n = h_n + p_{n-1}, \quad h_n = p_{n-1} + p_{n-2};$$

$$p_n = \frac{1}{2\sqrt{2}}\left[\left(1+\sqrt{2}\right)^{n+1} - \left(1-\sqrt{2}\right)^{n+1}\right], \quad h_n = \frac{1}{2}\left[\left(1+\sqrt{2}\right)^{n+1} + \left(1-\sqrt{2}\right)^{n+1}\right].$$

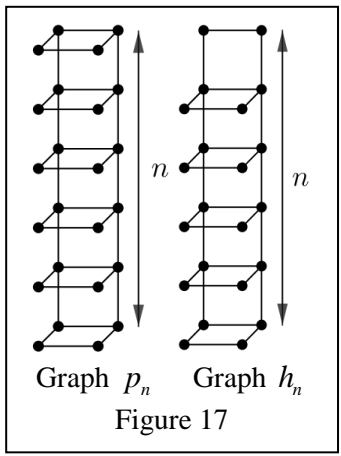

Graph $p_n$   Graph $h_n$
Figure 17

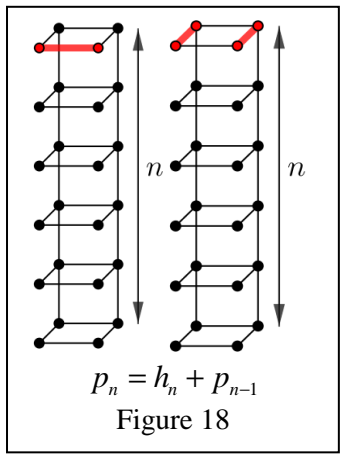

$p_n = h_n + p_{n-1}$
Figure 18

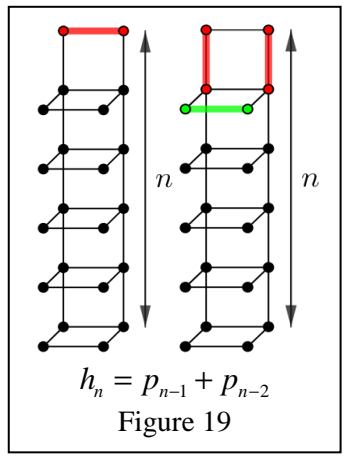

$h_n = p_{n-1} + p_{n-2}$
Figure 19

**Proof.** On figure 18 and figure 19 the relations between the domino tilings of graph $p_n$ and graph $h_n$ are specified. We have the system $p_n = h_n + p_{n-1}$ and $h_n = p_{n-1} + p_{n-2}$. Hence, we determine several dependencies:

$$p_n = \frac{h_n - h_{n-1}}{2}, \quad h_n = h_{n-1} + 2p_{n-2},$$

$$p_n = 2p_{n-1} + p_{n-2} \text{ и } h_n = 2h_{n-1} + h_{n-2}.$$

Furthermore, $p_1 = 2$, $p_2 = 5$, $p_3 = 12$, $h_1 = 1$, $h_2 = 3$, $h_3 = 7$, which means that $p_n$ are Pell's numbers, included in OEIS with number A000129, and the numbers $h_n$ are with number A001333.

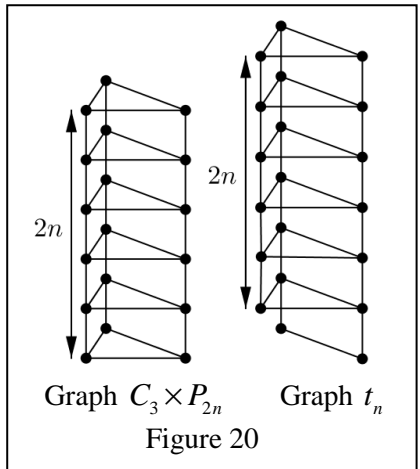

Graph $C_3 \times P_{2n}$   Graph $t_n$
Figure 20

### 5. Domino tilings of graph grids $C_3 \times P_{2n}$ and $K_{1,3} \times P_{2n}$

*Two more graph grids are associated with impressive dependencies for the domino tilings. These are graph $C_3 \times P_{2n}$ and graph $K_{1,3} \times P_{2n}$.*

**Theorem 8.** Let's consider graph $C_3 \times P_{2n}$ and graph $t_n$ as shown on figure 20. Let $T_n$ be the number of the domino tilings for graph $C_3 \times P_{2n}$, and $t_n$ be the number of the domino tilings of graph $t_n$, where $n$ is a positive integer. Then for $n \geq 1$:

$T_n = T_{n-1} + 3t_{n-1}$ ; $t_n = T_n + t_{n-1}$;

$T_n = 5T_{n-1} - T_{n-2}$ ; $t_n = 5t_{n-1} - t_{n-2}$;

$$T_n = \frac{1}{14}\left[(7+\sqrt{21})\left(\frac{5+\sqrt{21}}{2}\right)^n + (7-\sqrt{21})\left(\frac{5-\sqrt{21}}{2}\right)^n\right] \text{ for } n \geq 1;$$

$$t_n = \frac{\sqrt{21}}{21}\left[\left(\frac{5+\sqrt{21}}{2}\right)^{n+1} - \left(\frac{5-\sqrt{21}}{2}\right)^{n+1}\right] \text{ for } n \geq 0.$$

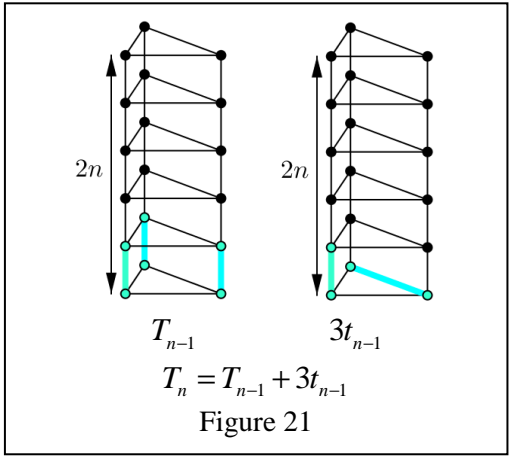

$T_{n-1}$   $3t_{n-1}$
$T_n = T_{n-1} + 3t_{n-1}$
Figure 21

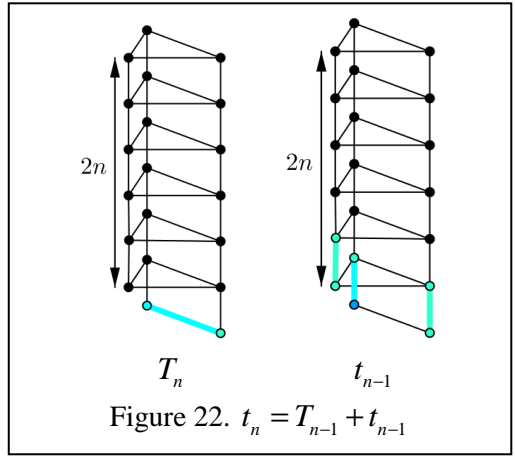

$T_n$   $t_{n-1}$
Figure 22. $t_n = T_{n-1} + t_{n-1}$

**Proof.** Let's start the counting from the "bottom" of the graph grid. The non-zero possibilities are illustrated on figure 21: three vertical edges, one vertical and one horizontal edge. We receive the first dependency $T_n = T_{n-1} + 3t_{n-1}$. We conclude that we have the recurrent equation $t_n = T_n + t_{n-1}$ (figure 22). From the obtained two recurrent equations for $T_n$ and $t_n$ it follows that $T_n = 5T_{n-1} - T_{n-2}$ and $t_n = 5t_{n-1} - t_{n-2}$. The sequence $\{T_n\}$ is linear-recurrent, and for this sequence we have $T_1 = 4$, $T_2 = 19$, the characteristic equation is $x^2 - 5x + 1 = 0$. Its roots are $x_1 = \frac{5+\sqrt{21}}{2}$ and $x_2 = \frac{5-\sqrt{21}}{2}$. Then

$$T_n = \frac{1}{14}\left[(7+\sqrt{21})\left(\frac{5+\sqrt{21}}{2}\right)^n + (7-\sqrt{21})\left(\frac{5-\sqrt{21}}{2}\right)^n\right], n \geq 1.$$

The sequence $\{t_n\}$ has the same characteristic equation and knowing that $t_0 = 1$, $t_1 = 5$, $t_2 = 24$ we receive for $n \geq 0$

$$t_n = \frac{\sqrt{21}}{21}\left[\left(\frac{5+\sqrt{21}}{2}\right)^{n+1} - \left(\frac{5-\sqrt{21}}{2}\right)^{n+1}\right].$$

*We will prove one unique coincidence for the number of the domino tilings of the two graphs $C_3 \times P_{2n}$ and $K_{1,3} \times P_{2n}$.*

**Theorem 9.** Let's consider graph $K_{1,3} \times P_{2n}$ and graph $q_n$ as shown on figure 23. Let $Q_n$ be the number of the domino tilings of graph $K_{1,3} \times P_{2n}$, and $q_n$ be the number of the domino tilings of graph $q_n$, where $n$ is a positive integer. Then for the number of the domino tilings of graph $C_3 \times P_{2n}$ and the number of the domino tilings of graph $K_{1,3} \times P_{2n}$ we have the equations: $Q_n = T_n$; $q_n = t_n$.

**Proof.** Let's start counting from the "bottom" of the graph grid $K_{1,3} \times P_{2n}$ (Another denotation is $S_3 \times P_{2n}$). On figure 24 are illustrated the non-zero possibilities: three vertical edges or one horizontal and two vertical edges. We obtain $Q_n = Q_{n-1} + 3q_{n-1}$. Now let's consider $q_n$ - using similar reasoning we conclude that we have the recurrent equation $q_n = Q_n + q_{n-1}$ (figure 25). From the obtained two recurrent equations for $Q_n$ and $q_n$ it follows that $Q_n = 5Q_{n-1} - Q_{n-2}$ and $q_n = 5q_{n-1} - q_{n-2}$.

The sequences $\{Q_n\}$ and $\{q_n\}$ are linear-recurrent and $Q_1 = 4$, $Q_2 = 19$, $q_1 = 5$, $q_2 = 24$. This means that $Q_n = T_n$ and $q_n = t_n$.

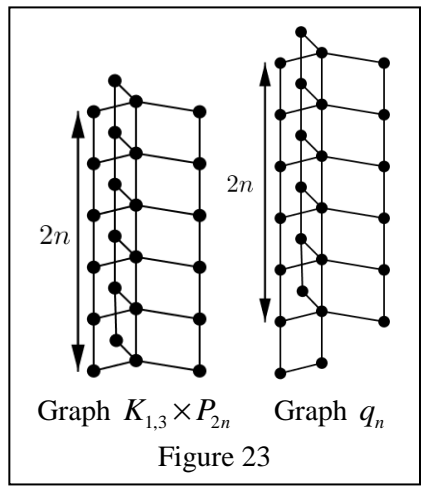

Graph $K_{1,3} \times P_{2n}$    Graph $q_n$

Figure 23

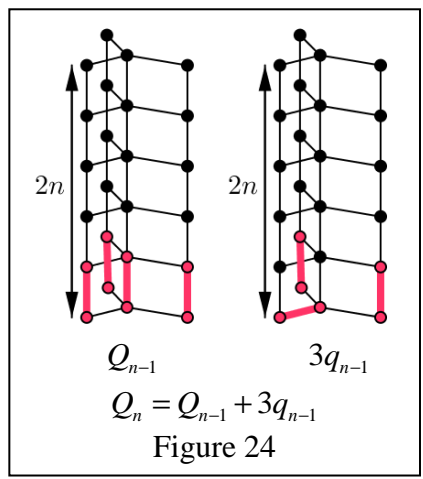

$Q_{n-1}$    $3q_{n-1}$

$Q_n = Q_{n-1} + 3q_{n-1}$

Figure 24

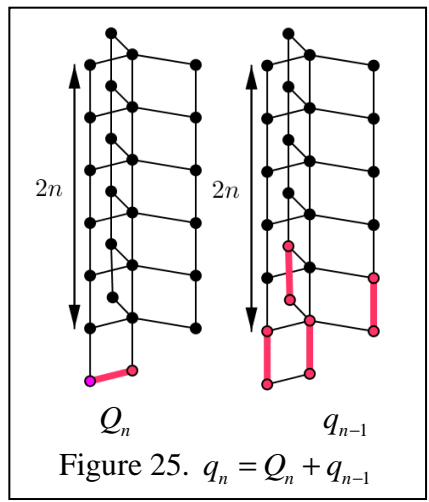

$Q_n$    $q_{n-1}$

Figure 25. $q_n = Q_n + q_{n-1}$

**Note.** The sequences $\{T_n\}$ and $\{t_n\}$, respectively $\{Q_n\}$ and $\{q_n\}$, are denoted in OEIS with numbers A004253 and A004254. For sequence A004253 it is noted that it represents the domino tilings of a graph, but for sequence A004254 it is not noted that it also represents the domino tilings of graph grids.

## 6. Domino tilings of graph grids $K_4 \times P_n$

**Theorem 10.** The domino tilings of graph $K_4 \times P_n$ (figure 26) are

$$M_n = \frac{1}{7}\left[\left(\frac{5+\sqrt{21}}{2}\right)^{n+1} + \left(\frac{5-\sqrt{21}}{2}\right)^{n+1} + 2(-1)^n\right].$$

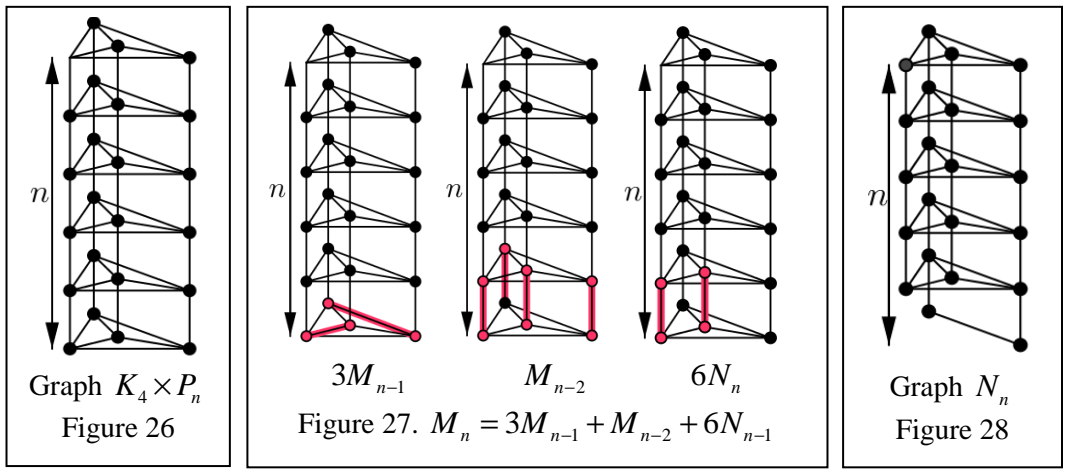

Graph $K_4 \times P_n$
Figure 26

$3M_{n-1}$  $M_{n-2}$  $6N_n$

Figure 27. $M_n = 3M_{n-1} + M_{n-2} + 6N_{n-1}$

Graph $N_n$
Figure 28

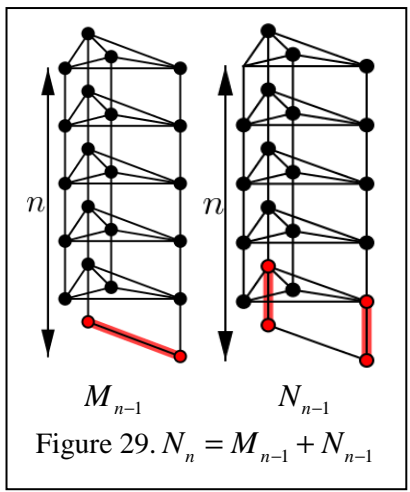

$M_{n-1}$  $N_{n-1}$

Figure 29. $N_n = M_{n-1} + N_{n-1}$

**Proof.** Let's start the counting from the "bottom" of the grid. On figure 27 the non-zero possibilities are illustrated: zero vertical edges, four vertical edges and two vertical edges. We note with $N_n$ the number of the domino tilings of graph $K_4 \times P_n$, from which two vertices and the corresponding edges are removed (figure 28). It follows that we receive the dependency $M_n = 3M_{n-1} + M_{n-2} + 6N_{n-1}$. We conclude that we have the recurrent equation $N_n = M_{n-1} + N_{n-1}$ (figure 29). We obtain a system of two recurrent equations $M_n = 3M_{n-1} + M_{n-2} + 6N_n$, $N_n = M_{n-1} + N_{n-1}$.

Hence, we have:
$M_n = 4M_{n-1} + 4M_{n-2} - M_{n-3}$;

$$N_n = 4N_{n-1} + 4N_{n-2} - N_{n-3}.$$

For sequence $\{M_n\}$ we have $M_1 = 3$, $M_2 = 16$, $M_3 = 75$, and for sequence $\{N_n\}$ - $N_1 = 1$, $N_2 = 4$, $N_3 = 20$. The characteristic equation of the two sequences is $x^3 - 4x^2 - 4x + 1 = 0$, its roots are $x_1 = \frac{5+\sqrt{21}}{2}$, $x_2 = \frac{5-\sqrt{21}}{2}$ and $x_3 = 1$. Then

$$M_n = \frac{1}{7}\left(\frac{5+\sqrt{21}}{2}\right)^{n+1} + \frac{1}{7}\left(\frac{5-\sqrt{21}}{2}\right)^{n+1} + \frac{2}{7}(-1)^n;$$

$$N_n = \frac{3+\sqrt{21}}{42}\left(\frac{5+\sqrt{21}}{2}\right)^n + \frac{3-\sqrt{21}}{42}\left(\frac{5-\sqrt{21}}{2}\right)^n + \frac{1}{7}(-1)^n.$$

*Let's consider a remarkable relation between the domino tilings of graph $K_4 \times P_n$ and graph $C_3 \times P_{2n}$, which is an interesting analogue of the relation between the domino tilings of graph $C_4 \times P_n$ and graph $P_3 \times P_{2n}$ from theorem 3.*

**Theorem 11.** Let $M_n$ be the number of the domino tilings of graph $K_4 \times P_n$, $T_n$ be the number of the domino tilings of graph $C_3 \times P_{2n}$, and $t_n$ be the number of the domino tilings of the graph, defined on figure 20. Then: $M_{2n} = T_n^2$; $M_{2n-1} = 3t_n^2$.

**Proof.** The statement is illustrated on table 6. We have

$$T_n^2 = \left\{\frac{1}{14}\left[(7+\sqrt{21})\left(\frac{5+\sqrt{21}}{2}\right)^n + (7-\sqrt{21})\left(\frac{5-\sqrt{21}}{2}\right)^n\right]\right\}^2 =$$

$$= \left(\frac{1}{14}\right)^2\left[(7+\sqrt{21})^2\left(\frac{5+\sqrt{21}}{2}\right)^{2n+2} + (7-\sqrt{21})^2\left(\frac{5-\sqrt{21}}{2}\right)^{2n+2}\right] =$$

$$= \frac{1}{7}\left[\left(\frac{5+\sqrt{21}}{2}\right)^{2n+1} + \left(\frac{5-\sqrt{21}}{2}\right)^{2n+1} + 2\right] = M_{2n}$$

$$3t_n^2 = 3 \cdot \left\{\frac{\sqrt{21}}{21}\left[\left(\frac{5+\sqrt{21}}{2}\right)^{n+1} - \left(\frac{5-\sqrt{21}}{2}\right)^{n+1}\right]\right\}^2 =$$

$$= 3 \cdot \frac{1}{21}\left[\left(\frac{5+\sqrt{21}}{2}\right)^{2n+2} - 2 + \left(\frac{5-\sqrt{21}}{2}\right)^{2n+2}\right] = M_{2n+1}.$$

| $n$ | $T_n$ | $t_n$ | $M_n$ |
|---|---|---|---|
| 1 | 4 | 1 | 3 |
| 2 | 19 | 5 | 16 |
| 3 | 91 | 24 | 75 |
| 4 | 436 | 115 | 361 |
| 5 | 2089 | 551 | 1728 |
| 6 | 10009 | 2640 | 8281 |
| 7 | 47956 | 12649 | 39675 |
| 8 | 229771 | 60605 | 190096 |

Table 6

| Sequences in theorems 8, 9 and 10 | | Sequences in OEIS | Number |
|---|---|---|---|
| $T_n$ $Q_n$ | 4, 19, 91, ... $T_n = 5T_{n-1} - T_{n-2}$ | 1, 4, 19, 91, 436, 2089, 10009, 47956, 229771, ... | **A004253** |
| $t_n$ $q_n$ | 1, 5, 24, 115, ... $t_n = 5t_{n-1} - t_{n-2}$ | 0, 1, 5, 24, 115, 551, 2640, 12649, 60605, 290376, ... | **A004254** |
| $M_n$ | 3, 16, 75, 361, ... $M_n = 4M_{n-1} + 4M_{n-2} - M_{n-3}$ | 3, 16, 75, 361, 1728, 8281, 39675, 190096, 910803, ... | **A003769** |
| $N_n$ | 1, 4, 20, 95,... $N_n = 4N_{n-1} + 4N_{n-2} - N_{n-3}$ | 1, 4, 20, 95, 456, 2184, 10465, 50140, ... | **A099025** |
| Table 7 | | | |

In table 7 information about the sequences from the previous theorems is given.

**Note**. The sequence $\{M_n\}$ is refernced in OEIS with number A003769, and it is noted that it illustrates the domino tilings of graph $K_4 \times P_n$, but the relation between the domino tilings of graph $K_4 \times P_n$ with graph $C_3 \times P_{2n}$ and graph $t_n$ is not noted.

We will also consider the relation of $N_n$ with two number sequences. This can be seen on table 8. The

| $k$ | $T_k$ | $t_k$ | $k$ | | $N_n$ | $n$ |
|---|---|---|---|---|---|---|
| | | 1 | 0 | | 1 | 1 |
| 1 | 4 | 1 | 0 | 4.1 | 4 | 2 |
| 1 | 4 | 5 | 1 | 4.5 | 20 | 3 |
| 2 | 19 | 5 | 1 | 19.5 | 95 | 4 |
| 2 | 19 | 24 | 2 | 19.24 | 456 | 5 |
| 3 | 91 | 24 | 2 | 91.24 | 2184 | 6 |
| 3 | 91 | 115 | 3 | 91.115 | 10465 | 7 |
| 4 | 436 | 115 | 3 | 436.115 | 50140 | 8 |
| 4 | 436 | 551 | 4 | 436.551 | 240236 | 9 |
| 5 | 2089 | 551 | 4 | 2089.551 | 1151039 | 10 |
| 5 | 2089 | 2640 | 5 | 2089.2640 | 5514960 | 11 |
| Table 8 | | | | | | |

numbers $N_n$ are obtained by multiplying the numbers $T_n$ and $t_n$ from theorem 8 in a defined order.

**Theorem 12.** We will define for $n \geq 1$ a sequence of numbers $w_n$ as follows: $w_{2k-1} = t_{k-1}$ and $w_{2k} = T_k$ for $k \geq 1$. Then $N_n = w_n \cdot w_{n+1}$ for $n \geq 2$.

**Proof.** In accordance with the definition for $w_n$, $N_n = T_k \cdot t_k$ when $n$ is an odd number, and $N_n = T_k \cdot t_{k-1}$ when $n$ is an even number. For the respective products we have

$$T_k \cdot t_{k-1} = \frac{1}{14}\left[\left(7+\sqrt{21}\right)\left(\frac{5+\sqrt{21}}{2}\right)^k + \left(7-\sqrt{21}\right)\left(\frac{5-\sqrt{21}}{2}\right)^k\right].$$

$$\cdot \frac{\sqrt{21}}{21}\left[\left(\frac{5+\sqrt{21}}{2}\right)^k - \left(\frac{5-\sqrt{21}}{2}\right)^k\right] =$$

$$= \frac{3+\sqrt{21}}{42}\left(\frac{5+\sqrt{21}}{2}\right)^{2k+1} + \frac{3-\sqrt{21}}{42}\left(\frac{5-\sqrt{21}}{2}\right)^{2k+1} - \frac{1}{7} = N_{2k}.$$

$$T_k \cdot t_k = \frac{1}{14}\left[\left(7+\sqrt{21}\right)\left(\frac{5+\sqrt{21}}{2}\right)^k + \left(7-\sqrt{21}\right)\left(\frac{5-\sqrt{21}}{2}\right)^k\right].$$

$$\cdot \frac{\sqrt{21}}{21}\left[\left(\frac{5+\sqrt{21}}{2}\right)^{k+1} - \left(\frac{5-\sqrt{21}}{2}\right)^{k+1}\right] =$$

$$= \frac{3+\sqrt{21}}{42}\left(\frac{5+\sqrt{21}}{2}\right)^{2k+1} + \frac{3-\sqrt{21}}{42}\left(\frac{5-\sqrt{21}}{2}\right)^{2k+1} + \frac{1}{7} = N_{2k+1}.$$

**Note.** The sequence $\{N_n\}$ is referenced in OEIS with number A099025, but it is not noted that this is also the number of the domino tilings of a graph, which is obtained when two vertices are removed from graph $K_4 \times P_n$, as shown on figure 28.

**26.12.2018**

Valcho Milchev, teacher
Kardzhali, Bulgaria
e-mail: milchev.vi@gmail.com,   milchev_v@abv.bg